\documentclass{article}
\usepackage{indentfirst,latexsym,bm,amscd,amsmath,amsthm,amsfonts,array,amssymb,eucal}
\begin{document}
\title{The Monsky-Washnitzer cohomology and the\\ de Rham cohomology}
\author{BinYong Hsie\\LMAM, Department of Mathematics, Peking university\\Bejing, 100871, P.R.China}
\maketitle
\abstract{The author constructs a theory of dagger formal schemes
over $R$ and then defines the de Rham cohomology for flat dagger
formal schemes $X$ with integral and regular reductions $\bar{X}$
which generalizes the Monsky-Washnitzer cohomology. Finally the
author gets Lefschetz' fixed pointed formula for $X$ with certain
conditions.}
\section{Introduction}
Let $K$ be a finite extension of $\mathbb{Q}_{p}$ with $R$ its
ring of integers and $k=\mathbb{F}_{q}$ its residue field. Let
$\pi$ be a uniformizer of $R$.

Let $\bar{X}$ be a regular algebraic variety over $k$. Define the
zeta-function $Z(\bar{X}|k,t)$ of $\bar{X}$ by
$Z(\bar{X}|k,t)=\mathrm{exp}(\sum\limits_{s\geq 1}\frac{N_{s}}{s}
t^{s})$ with $N_{s}$ the number of $\mathbb{F}_{q^{s}}$-points of
$\bar{X}$.

Weil's conjecture says that $Z(\bar{X}|k,t)$ is a rational
function. To prove it, one tries to find suitable cohomology such
that the Lefschetz' fixed points formula holds. When $\bar{X}$ is
an affine integral and regular variety over $k$, the
Monsky-Washnitzer cohomology is such a cohomology.

Let $\bar{X}=\mathrm{Spec}(\bar{A})$ with $\bar{A }$ integral and
regular over $k$. Let $A$ be a flat w.c.f.g. algebra over
$R=\mathrm{W}(k)$ which is a lift of $\bar{A}$. Note that every
flat lift of $\bar{A}$ is $R$-isomorphic to $A$, and that
$\Omega^{1}(A)=\Omega^{1}(A/R)$ is a projective $A$-module. One
can define the de Rham complex $\Omega^{\cdot}(A)$
$$0\rightarrow \Omega^{0}(A)\xrightarrow[d^{0}]{} \Omega^{1}(A)\xrightarrow[d^{1}]{} \Omega^{2}(A)\rightarrow\cdot\cdot\cdot$$
with $\Omega^{i}(A)=\bigwedge^{i}\Omega^{1}(A)$.
$H^{i}_{\small\mathrm{MW}}(\bar{X},
K):=H^{i}(\Omega(A))\otimes_{R}K$ is the definition of the
Monsky-Washnitzer cohomology. Note that
$H^{i}(\Omega(A))\otimes_{R}K=H^{i}(\Omega(A)\otimes_{R}K)$. If
$\mathrm{F}$ is a lift of the Frobenius map $x\mapsto x^{q}$ over
$\bar{A}$ to $A$, then $\mathrm{F}^{*}$ induces an isomorphism
over $H^{i}_{\small\mathrm{MW}}(\bar{X}, K)$ which does not depend
on the choice of lift $\mathrm{F}$. Moreover,
$(\mathrm{F}^{*})^{-1}$ is a nuclear operator, and satisfies
$$
N_{s}=\sum (-1)^{i} \mathrm{tr}
((q^{n}(\mathrm{F}^{*})^{-1})^{s}|H^{i}_{\small\mathrm{MW}}(\bar{X},
K)).
$$

We find a generalization of the Monsky-Washnitzer cohomology. We
construct a theory of dagger formal scheme over $R$ and define the
de Rham cohomology for flat dagger formal schemes $X$ with
integral and regular reductions.

An affine dagger formal scheme is a pair
$(\mathrm{Spec}(A\otimes_{R}k),\mathcal{O}^{\dagger})$, where $A$
is a w.c.f.g. algebra over $R$, and $\mathcal{O}^{\dagger}$ is a
sheaf over $\bar{X}=\mathrm{Spec}(A\otimes_{R}k)$ such that
$\mathcal{O}^{\dagger}(\bar{X}_{\bar{f}})=A\langle
f^{-1}\rangle^{\dagger}$. Let $\mathrm{Spf}^{\dagger}(A)$ denote
the pair $(\mathrm{Spec(A\otimes_{R}k)},\mathcal{O}^{\dagger})$.
$\mathcal{O}^{\dagger}$ is called the structure sheaf of
$\mathrm{Spf}^{\dagger}(A)$. One can show
$\mathrm{Spf}^{\dagger}(A)$ is a locally ringed space. A dagger
formal scheme is a locally ringed space
$(X,\mathcal{O}^{\dagger}_{X})$ in which every point has an open
neighborhood $U$, such that $(U,\mathcal{O}^{\dagger}_{X}|U)$ is
an affine dagger formal scheme.

When $X$ is flat, separated and Noetherian over $R$ with $\bar{X}
$ integral and regular, the sheaf $\Omega^{1}$ over $X$ is locally
free. So one can define the de Rham complex $\Omega^{\cdot}$ and
the de Rham cohomology $H^{i}_{dR}(X;
K):=\mathbb{H}^{i}(X,(\Omega^{\cdot}\otimes_{R}K))$. When
$X=\mathrm{Spf}^{\dagger}(A)$, we have
$H^{i}_{dR}(X,K)=H^{i}_{\small\mathrm{MW}}(\bar{X},K)$. When $R$
is $\mathrm{W}(k)$, and there is a lift $\mathrm{F}$ of the
Frobenius of $\bar{X}$ to $X$, we have the following theorem.

{\bf Theorem.} {\itshape Let $X$ be a flat separated and
Noetherian dagger formal scheme with $\bar{X}=X\otimes_{R}k$
regular and integral of dimension $n$.  Then $\mathrm{F}^{*}$ is
an isomorphism, $(\mathrm{F}^{*})^{-1}$ is nuclear over
$H_{dR}^{i}(X,K)$. And we have the following formula} $$
N_{s}(\bar{X})=\sum (-1)^{i} \mathrm{tr}
((q^{n}(\mathrm{F}^{*})^{-1})^{s}|H^{i}_{dR}(X, K)).
$$

\section{Weakly complete finitely generated algebra}
\subsection{Definition of weakly complete finitely generated algebra}
For a nonnegative integer $n$, let us define
$$\mathbb{T}_{n}:=\{\sum_{v \in \mathbb{N}^{n}}c_{v}\xi^{v}\in R[[\xi_{1},...,\xi_{n}]]: |c_{v}|\rightarrow 0\},$$
and
$$\mathbb{T}^{\dagger}_{n}:=\{\sum_{v \in \mathbb{N}^{n}}c_{v}\xi^{v}\in R[[\xi_{1},...,\xi_{n}]]: \exists \varepsilon> 0, |c_{v}|p^{\varepsilon |v|}\rightarrow 0\}.$$
It is well know that $\mathbb{T}_{n}$ is complete with respect to
the Gauss norm. There is also a Gauss norm over
$\mathbb{T}^{\dagger}_{n}$. But $\mathbb{T}^{\dagger}_{n}$ is not
complete with respect to this norm.

{\bf Proposition 1.}(\cite{Put}) $\mathbb{T}^{\dagger}_{n}$
{\itshape
satisfies Weierstrass preparation and division.} 

As a consequence, $\mathbb{T}^{\dagger}_{n}$ is noetherian and
flat over $R[\xi_{1},...,\xi_{n}]$.

{\bf Definition.} {\itshape A weakly complete finitely generated
(w.c.f.g.) algebra $A$ over $R$ is a homomorphic image of some
$\mathbb{T}^{\dagger}_{n}$. While, a complete finitely generated
(c.f.g.) algebra $\breve{A}$ over $R$ is a homomorphic image of
some $\mathbb{T}_{n}$.}

\subsection{Faithful flat of $\breve{A}$ over $A$}
Let $A$ be a neotherian ring with $I$ an ideal of $A$. Then the
$I$-adic completion $\breve{A}$ of $A$ is flat over $A$.

In this subsection, we need to consider when $\breve{A}$ is a
faithful flat over $A$. For this, we have the following theorem.

{\bf Theorem 1.}(\cite{Matsumura}) {\itshape Let $A$ be a
neotherian ring with an adic topology, and let $I$ be an ideal of
definition. Then the following is equivalent.

 (1). $A$ is a Zariski ring, i.e., every ideal is closed in it.

 (2). $I\subseteq rad(A)$.

 (3). Every finite $A$-module $M$ is separated in the $I$-adic
 topology.

 (4). In every finite $A$-module $M$, every submodule is closed in
 the $I$-adic topology.

 (5). The completion $\breve{A}$ of $A$ is faithful flat over $A$.}

We have the following lemma.

{\bf Lemma 1.} {\itshape A w.c.f.g. algebra $A$ over $R$ is a
Zariski ring.}

{\bf Proof.} It is enough to show that the lemma holds when
$A=\mathbb{T}^{\dagger}_{n}$.

For every $f\in \mathbb{T}^{\dagger}_{n}$, let us show $1-\pi f$
is invertible. In fact, we only need to show
$\sum_{n=0}^{\infty}(\pi f)^{n}\in \mathbb{T}^{\dagger}_{n}$. Let
$f=\sum_{\alpha\in \mathbb{N}^{n}}a_{\alpha}\xi^{\alpha}$. Since
$f\in \mathbb{T}^{\dagger}_{n}$, there are two positive numbers
$M$ and $\varepsilon$ such that $|a_{\alpha}|< Mp^{-\varepsilon
|\alpha|}$. Therefore, there is a positive integer $N$ such that
when $|\alpha|\geq N$,
$|a_{\alpha}|<p^{-\frac{\varepsilon}{2}|\alpha|}$.

Let $\varepsilon'=
\mathrm{min}(\frac{\varepsilon}{2},\frac{-\mathrm{log}_{p}|\pi|}{N})$
and $\pi f=\sum_{\alpha\in
\mathbb{N}^{n}}a'_{\alpha}\xi^{\alpha}$. Then
\begin{equation}|a'_{\alpha}|<p^{-\varepsilon'|\alpha|}.\label{estimite}\end{equation}
It is easy to show $\sum_{n=0}^{\infty}(\pi f)^{n}$ converges in
$R[[\xi_{1},...,\xi_{n}]]$. Let $\sum_{n=1}^{\infty}(\pi
f)^{n}=\sum_{\alpha\in \mathbb{N}^{n}}b_{\alpha}\xi^{\alpha}$.
Then from fomula (\ref{estimite}), we know
$|b_{\alpha}|<p^{-\varepsilon'|\alpha|}$. Therefore,
$\sum_{n=0}^{\infty}(\pi f)^{n}\in \mathbb{T}^{\dagger}_{n}$ as
desired. \qed

Now, the following proposition is easily deduced from Lemma 1 and
Theorem 1.

{\bf Proposition 2.} {\itshape Let $A$ be a w.c.f.g. algebra with
$\breve{A}$ its completion according to the $\pi$-adic topology.
Then $\breve{A}$ is a faithful flat $A$-algebra and $A\rightarrow
\breve{A}$ is universal injective.}

Let $A$ be a w.c.f.g. algebra with $f \in A$. Let us define
$A\langle f^{-1}\rangle^{\dagger}$ to be $A\langle
\xi\rangle^{\dagger}/(f\xi -1)$. Then $A\langle
f^{-1}\rangle^{\dagger}$ is a w.c.f.g. algebra. If $f,g$ are two
elements in $A$ such that $f-g \in \pi A$, then $A\langle
f^{-1}\rangle^{\dagger}$ is canonically isomorphic to $A\langle
g^{-1}\rangle^{\dagger}$. We have the following proposition.

{\bf Proposition 3.} {\itshape $A\langle f^{-1}\rangle^{\dagger}$
is a flat $A$-algebra.}

 {\bf Proof.} It is well know that $A_{f^{-1}}$ is a
flat $A$-algebra. $A_{f^{-1}}$ and $A\langle
f^{-1}\rangle^{\dagger}$ have the same $\pi$-adic completion $B$.
$B$ is a flat $A_{f^{-1}}$-algebra and is a faithful flat
$A\langle f^{-1}\rangle^{\dagger}$-algebra according to
Proposition 2. Therefore, $A\langle f^{-1}\rangle^{\dagger}$ is a
flat $A$-algebra. \qed

From now on, we always write $\breve{A}$ for the $\pi$-adic
completion of a w.c.f.g. $R$-algebra $A$.

Let $A_{\lambda}=A\otimes_{R}(R/\pi^{\lambda}R)$.

{\bf Proposition 4.} {\itshape Let $\varphi:\: A\rightarrow B$ be
a morphism of w.c.f.g. $R$-algebras. And let $M$ be a finite
$B$-module. Then $M$ is (faithful) flat over $A$ if and only if
$M\otimes_{R}R_{\lambda}$ is (faithful) flat over
$A\otimes_{R}R_{\lambda}$ for all $\lambda\in \mathbb{N}$.}

Its proof is similar to that in \cite{Bosch I} Lemma 1.6.

Let $A$ be a flat w.c.f.g. algebra over $R$. Let $f$ be an element
of $A$ with $f_{1},...,f_{r}\in (f)$. Assume that
$$\overline{f}\in \sqrt{(\overline{f_{1}},...,\overline{f_{r}})}\eqno(*).$$
Then we have the following corollary of Proposition 4.

{\bf Corollary 1.} {\itshape $\prod\limits_{i}A\langle
f_{i}^{-1}\rangle^{\dagger}$ is faithful flat over $A\langle
f^{-1}\rangle^{\dagger}$.}

If $A$ is a flat w.c.f.g. $R$-algebra, then both $A\rightarrow
A\otimes_{R} K$ and $\breve{A}\rightarrow \breve{A}\otimes_{R} K$
are injective. Moreover, we have the following lemma.

{\bf Lemma 2.} {\itshape In $\breve{A}\otimes_{R} K$, the
intersection of $A\otimes_{R} K$ and $\breve{A}$ is $A$.}

{\bf Proof.} Since $A$ is a Zariski ring, every ideal $I$ of $A$
is closed in $A$. So $I\breve{A}\cap A=I$. See \cite{Zariski} for
more details.

For any $\frac{a}{\pi^{n}}\in A\otimes K$ with $a\in A$, if
$\frac{a}{\pi^{n}}\in \breve{A}$, then $a\breve{A}\subseteq \pi
^{n}\breve{A}$. We get $aA\subseteq \pi^{n} A$. Then
$\frac{a}{\pi^{n}}\in A$ as desired. \qed

\subsection{Dagger rigid geometry}
We can construct a theory of rigid geometry using dagger affinoid
algebras over $K$ instead of affinoid algebras.

{\bf Definition.} {\itshape A dagger affinoid algebra over $K$ is
a homomorphic image of some
$\mathbb{T}^{\dagger}_{n}\otimes_{R}K$.}

There is a one to one correspondence between prime ideals of
$\mathbb{T}^{\dagger}_{n}\otimes_{R}K$ and prime ideals of
$\mathbb{T}^{\dagger}_{n}$ that do not contain $\pi$. The same
assertion holds if we use $\mathbb{T}^{\dagger}_{n}\otimes_{R}K$
and $\mathbb{T}^{\dagger}_{n}$ instead of
$\mathbb{T}_{n}\otimes_{R}K$ and $\mathbb{T}_{n}$, respectively.
Since $\mathbb{T}_{n}$ is faithful flat over
$\mathbb{T}^{\dagger}_{n}$, every (maximal) prime ideal of
$\mathbb{T}^{\dagger}_{n}\otimes_{R}K$ is the restriction of a
prime (maximal) ideal of $\mathbb{T}_{n}\otimes_{R}K$. Moreover,
if $I$ is an ideal of $\mathbb{T}^{\dagger}_{n}\otimes_{R}K$, then
$I=I(\mathbb{T}_{n}\otimes_{R}K)\cap(\mathbb{T}^{\dagger}_{n}\otimes_{R}K)$.

Let $\mathcal{M}^{\mathrm{v}}$ be a maximal ideal of
$\mathbb{T}_{n}\otimes_{R}K$, then
$\mathbb{T}_{n}\otimes_{R}K/\mathcal{M}^{\mathrm{v}}$ is an
algebraic extension of $K$. So the restriction of
$\mathcal{M}^{\mathrm{v}}$ to
$\mathbb{T}^{\dagger}_{n}\otimes_{R}K$ is again a maximal ideal.
Since $\mathbb{T}^{\dagger}_{n}\otimes_{R}K$ is dense in
$\mathbb{T}_{n}\otimes_{R}K$, two different maximal ideals have
different restrictions over
$\mathbb{T}^{\dagger}_{n}\otimes_{R}K$. Therefore,
$\mathbb{T}^{\dagger}_{n}\otimes_{R}K$ and
$\mathbb{T}_{n}\otimes_{R}K$ have the same maximal spectral. Let
$B=\mathbb{T}^{\dagger}_{n}\otimes_{R}K/I$, and let
$I^{\mathrm{v}}$ be $I(\mathbb{T}_{n}\otimes_{R}K)$. Then
$I=I^{\mathrm{v}}\cap (\mathbb{T}^{\dagger}_{n}\otimes_{R}K)$. Let
$B^{\mathrm{v}}=\mathbb{T}_{n}\otimes_{R}K/I^{\mathrm{v}}$, then
$B\rightarrow B^{\mathrm{v}}$ is injective. And there is a one to
one correspondence between $\mathrm{Spm}(B)$ and
$\mathrm{Spm}(B^{\mathrm{v}})$. Moveover $B^{\mathrm{v}}$ is the
completion of $B$ according to the maximal spectral.

We can define Weierstrass domains, Laurent domains and rational
domains for $X=\mathrm{Spm}(B)$, and then define Grothendieck
topology for $X$. We use $X^{\mathrm{v}}$ to denote
$\mathrm{Spm}(B^{\mathrm{v}})$.

{\bf Definition.} {\itshape

$(\mathrm{i})$. A subset in $X$ of type
$$X(f_{1},...,f_{r}):=\{x\in X: |f_{i}(x)|\leq 1\}$$
for $f_{1},...,f_{r}\in B$ is called a Weierstrass domain in $X$.

$(\mathrm{ii})$. A subset in $X$ of type
$$X(f_{1},...,f_{r}, g_{1}^{-1},...,g_{s}^{-1}):=\{x\in X: |f_{i}(x)|\leq 1; |g_{j}(x)|\geq 1\}$$
for $f_{1},...,f_{r}, g_{1},...,g_{s}\in B$ is called a Laurent
domain in $X$.

$(\mathrm{iii})$. A subset in $X$ of type
$$X(\frac{f_{1}}{f_{0}}, ..., \frac{f_{r}}{f_{0}} ):=\{ x\in X: |f_{i}(x)| \leq |f_{0}(x)| \}$$
for $f_{0},...,f_{r}\in B$ without common zero is called a
rational domain in $X$.}

If $f_{1}',...,f_{r}', g_{1}',...,g_{s}'$ are all nonzero elements
in $ B^{\mathrm{v}}$, choose $f_{1},...,f_{r},
\\g_{1},...,g_{s}\in B$ such that $\parallel
f_{i}-f_{i}'\parallel_{max}\leq 1$, $\parallel
g_{i}-g_{i}'\parallel_{max}\leq 1$. Then
$$X(f_{1},...,f_{r})=X^{\mathrm{v}}(f_{1}',...,f_{r}'),$$
and
$$X(f_{1},...,f_{r}, g_{1}^{-1},...,g_{s}^{-1})=X^{\mathrm{v}}(f_{1}',...,f_{r}', g_{1}'^{-1},...,g_{s}'^{-1}).$$
Assume $f_{0}',...,f_{r}'\in B^{\mathrm{v}}$ have no common zero.
Let $U_{i}^{\mathrm{v}}$ denote
$X^{\mathrm{v}}(\frac{f_{0}'}{f_{i}'}, ...,
\frac{f_{r}'}{f_{i}'})$. Then there is a positive number
$\varepsilon$ such that $|f_{i}' (x)|
> \varepsilon$ for any $x \in U_{i}^{\mathrm{v}}$. Choose $f_{0},...,f_{r}\in
B$ such that $\parallel f_{i}-f_{i}'\parallel_{max}<\varepsilon$,
then
$$U_{i}^{\mathrm{v}}=X(\frac{f_{1}}{f_{0}},...,\frac{f_{r}}{f_{0}}).$$
Therefore, we can define the Grothendieck topology for
$\mathrm{Spm}(B)$ such that it has the same Grothendieck topology
as $\mathrm{Spm}(B^{\mathrm{v}})$.

We can get Tate's acyclic theorem for dagger affinoid algebras in
the same way as in \cite{Bosch} or as in \cite{BGR}.

Let $A$ be a flat w.c.f.g. algebra over $R$. Let $B=A\otimes_{R}K$
and $X=\mathrm{Spm}(B)$. Let $f$ be an element of $A$ with
$f_{1},...,f_{r}\in (f)$. If we have $\overline{f}\in
\sqrt{(\overline{f_{1}},...,\overline{f_{r}})}$ in
$A\otimes_{R}k$, then $\{X(f_{i}^{-1})\}_{i=1}^{r}$ is a finite
affine covering of $X(f^{-1})$. Therefore as a consequence of
Tate's acyclic theorem, we get the following exact sequence
\begin{equation}
0\rightarrow A\langle f^{-1}\rangle^{\dagger} \otimes_{R} K
\rightarrow\prod_{i}A\langle
f_{i}^{-1}\rangle^{\dagger}\otimes_{R} K\rightrightarrows
\prod_{i,j}A\langle (f_{i}f_{j})^{-1}\rangle^{\dagger}\otimes_{R}
K. \label{1}
\end{equation}

{\bf Proposition 5.} {\itshape For any $A$-module $M$, the
following sequence is exact}
\begin{equation}
0\rightarrow A\langle
f^{-1}\rangle^{\dagger}\otimes_{A}M\rightarrow
\prod\limits_{i}A\langle
f_{i}^{-1}\rangle^{\dagger}\otimes_{A}M\rightrightarrows\prod\limits_{i,j}A\langle
(f_{i}f_{j})^{-1}\rangle^{\dagger}\otimes_{A}M.
\end{equation}
{\bf Proof.} (i). We have the following commutative diagram
($\otimes=\otimes_{A}$ in this proof)
\[ {\scriptsize
\begin{CD}
0 @>>> A\langle f^{-1}\rangle^{\dagger}\otimes M   @>>>
\prod\limits_{i}A\langle f_{i}^{-1}\rangle^{\dagger}\otimes M @>>>
\prod\limits_{i,j}A\langle
(f_{i}f_{j})^{-1}\rangle^{\dagger}\otimes M\\
& & @VVV @VVV @VVV \\
0 @>>> (A\langle f^{-1}\rangle^{\dagger})^{\mathrm{v}}\otimes M
@>>> \prod\limits_{i}(A\langle
f_{i}^{-1}\rangle^{\dagger})^{\mathrm{v}}\otimes M @>>>
\prod\limits_{i,j}(A\langle
(f_{i}f_{j})^{-1}\rangle^{\dagger})^{\mathrm{v}}\otimes M.
\end{CD}}
\]
From the theory of formal schemes, we know the second line is
exact. From Proposition 2, we see all vertical maps are injective.
From Corollary 1, we get that $0\rightarrow A\langle
f^{-1}\rangle^{\dagger}\otimes M\rightarrow
\prod\limits_{i}A\langle f_{i}^{-1}\rangle^{\dagger}\otimes M$ is
exact. Therefore, to show the first line is exact, it is enough to
show that $(\prod\limits_{i}A\langle
f_{i}^{-1}\rangle^{\dagger}\otimes M)\cap ((A\langle
f^{-1}\rangle^{\dagger})^{\mathrm{v}}\otimes M)=A\langle
f^{-1}\rangle^{\dagger}\otimes M$.

(ii). From formula (\ref{1}), the second line of the above
commutative diagram (with $M=A$) and Lemma 2, we see when $M=A$,
the first line of the above commutative diagram is also exact.
Then $(\prod\limits_{i}A\langle f_{i}^{-1}\rangle^{\dagger})\cap
(A\langle f^{-1}\rangle^{\dagger})^{\mathrm{v}}=A\langle
f^{-1}\rangle^{\dagger}$.

(iii). Let $A'$ be the sum of $(A\langle
f^{-1}\rangle^{\dagger})^{\mathrm{v}}$ and
$\prod\limits_{i}A\langle f_{i}^{-1}\rangle^{\dagger}$ in
$\prod\limits_{i}(A\langle
f_{i}^{-1}\rangle^{\dagger})^{\mathrm{v}}$. Then $A'$ is a Zariski
ring, since it is a quotient of the Zariski ring $(A\langle
f^{-1}\rangle^{\dagger})^{\mathrm{v}}\times
\prod\limits_{i}A\langle f_{i}^{-1}\rangle^{\dagger}$. $A'$ is
dense in $\prod\limits_{i}(A\langle
f_{i}^{-1}\rangle^{\dagger})^{\mathrm{v}}$, so
$(A')^{\mathrm{v}}=\prod\limits_{i}(A\langle
f_{i}^{-1}\rangle^{\dagger})^{\mathrm{v}}$. Therefore
$\prod\limits_{i}(A\langle
f_{i}^{-1}\rangle^{\dagger})^{\mathrm{v}}$ is faithful flat over
$A'$ and $A'\rightarrow \prod\limits_{i}(A\langle
f_{i}^{-1}\rangle^{\dagger})^{\mathrm{v}}$ is universal injective.
We see that $A'$ is a flat $A$-algebra.

Then both \[{\tiny 0\rightarrow A\langle
f^{-1}\rangle^{\dagger}\otimes M\rightarrow(\oplus_{i}A\langle
f_{i}^{-1}\rangle^{\dagger}\otimes M)\oplus (A\langle
f^{-1}\rangle^{\dagger})^{\mathrm{v}}\otimes M \xrightarrow{\mbox{
p }} A'\otimes M\rightarrow 0 }\] with $p: (a,b)\mapsto a-b$, and
$$ 0\rightarrow A'\otimes M\rightarrow \prod\limits_{i}(A\langle
f_{i}^{-1}\rangle^{\dagger})^{\mathrm{v}}\otimes M $$ are exact.
So we get the following exact sequence
$${\tiny 0\rightarrow A\langle f^{-1}\rangle^{\dagger}\otimes M\rightarrow(\oplus_{i}A\langle f_{i}^{-1}\rangle^{\dagger}\otimes M)\oplus (A\langle f^{-1}\rangle^{\dagger})^{\mathrm{v}}\otimes M\rightarrow
\prod\limits_{i}(A\langle
f_{i}^{-1}\rangle^{\dagger})^{\mathrm{v}}\otimes M. }$$ Then we
obtain $\prod\limits_{i}A\langle
f_{i}^{-1}\rangle^{\dagger}\otimes M\cap (A\langle
f^{-1}\rangle^{\dagger})^{\mathrm{v}}\otimes M=A\langle
f^{-1}\rangle^{\dagger}\otimes M$ as desired. \qed

\section{Dagger formal schemes}
\subsection{Theory of $\breve{\mathrm{C}}$ech cohomology}
In this subsection, we recall the theory of
$\breve{\mathrm{C}}$ech cohomology. We follow \cite{Tamme}
closely.

Let $X$ be a (compact) topology space. Let $\mathcal{S}$ denote
the category of sheaves over $X$, and let $\mathcal{P}$ denote the
category of persheaves over $X$. Let $i: \mathcal{S}\rightarrow
\mathcal{P}$ denote the natural inclusion functor.

Let us consider the right derived functors of $i$, and define
$$\mathcal{H}^{q}(\cdot):=R^{q}i().$$

{\bf Proposition.} {\itshape For any abelian sheaf $F$ over $X$,
and any open subset $U$ of $X$, we have}
$$\mathcal{H}^{q}(F)(U)=H^{q}(U,F).$$

Let $\{U_{i}\rightarrow U\}_{i\in I}$ be a covering. For any
presheaf $F$ over $X$, we define a complex $\{\breve{C}^{q}\}$ by
$$\breve{C}^{q}(\{U_{i}\rightarrow U\}, F):=\prod_{i_{0}<...<i_{q}}F(U_{i_{0}}\cap\cdot\cdot\cdot\cap U_{i_{q}})$$
with $d^{q}: \breve{C}^{q}(\{U_{i}\rightarrow U\}, F)\rightarrow
\breve{C}^{q+1}(\{U_{i}\rightarrow U\}, F)$ defined by
$$(d^{q}s)_{i_{0}\cdot\cdot\cdot i_{q+1}}=\sum_{v=0}^{q+1}(-1)^{v}s_{i_{0}\cdot\cdot\cdot\hat{i_{v}}\cdot\cdot\cdot i_{q+1}}.$$

Now, we define 
$$\breve{H}^{q}(\{U_{i}\rightarrow U\}, F):=\mathrm{Ker}(d^{q})/\mathrm{Im}(d^{q-1}).$$

{\bf Proposition.} {\itshape There is a canonical isomorphism
between $\breve{H}^{q}(\{U_{i}\rightarrow U\}, F)$ and
$R^{q}\breve{H}^{0}(\{U_{i}\rightarrow U\}, F)$.}

We have the following theorem.

{\bf Theorem.} (Spectral sequence for $\breve{\mathrm{C}}$ech
cohomology)

{\itshape Let $\{U_{i}\rightarrow U\}$ be a covering. For each
$F\in\mathcal{S}$, there is a spectral sequence
$$E_{2}^{pq}=\breve{H}^{p}(\{U_{i}\rightarrow U\},\mathcal{H}^{q}(F))\Rightarrow E^{p+q}=H^{p+q}(U,F)$$
which is functorial in $F$.}

Let $\mathcal{U}$ be a family of open subsets of $X$ such that
\\(i) The intersection of two open subsets in $\mathcal{U}$ is in
$\mathcal{U}$ again, \\(ii) Each finite covering of an open subset
of $X$ has a refinement consisting of sets in $\mathcal{U}$.

Let $\{U_{i}\rightarrow U\}$be a covering. If $U\in\mathcal{U}$
and $U_{i}\in\mathcal{U}$, we call it a covering in $\mathcal{U}$.

We define $\mathcal{U}$-presheaves for $X$.

{\bf Definition.} {\itshape A $\mathcal{U}$-presheaf over $X$, is
a collection of abelian groups $F(U)$ for $U\in \mathcal{U}$, and
a collection of restriction morphisms $res_{U}^{V}:F(V)\rightarrow
F(U)$ for each $U\subset V$ such that
$res_{U}^{V}res_{V}^{W}=res_{U}^{W}$ for $U\subset V\subset W$ and
$res_{U}^{U}=id$. Let $\mathcal{P}_{\mathcal{U}}$ denote the
category of $\mathcal{U}$-presheaves.}

Let $j_{\mathcal{U}}$ be the natural functor from $\mathcal{P}$ to
$\mathcal{P}_{\mathcal{U}}$, then $j_{\mathcal{U}}$ is exact. Let
$i_{\mathcal{U}}=j_{\mathcal{U}}\circ i$.

{\bf Definition.} {\itshape An abelian sheaf $F$ over $X$ is
called $\mathcal{U}$-flabby if $\breve{H}^{q}(\{U_{i}\rightarrow
U\},i_{\mathcal{U}}(F))=0$ for $q> 0$ and each covering in
$\mathcal{U}$.}

{\bf Proposition 6.} {\itshape $\mathrm{i)}$. Let $0\rightarrow
F'\rightarrow F\rightarrow F''$ be an exact sequence in
$\mathcal{S}$. If $F'$ is $\mathcal{U}$-flabby, the sequence is
exact in $\mathcal{P}_{\mathcal{U}}$ as well.

$\mathrm{ii)}$. Let $0\rightarrow F'\rightarrow F\rightarrow F''$
be an exact sequence in $\mathcal{S}$. If $F'$  and $F$ are
$\mathcal{U}$-flabby, so is $F''$.

$\mathrm{iii)}$. If the direct sum $F\oplus G$ of abelian sheaves
is $\mathcal{U}$-flabby, so is $F$.

$\mathrm{iv)}$. Injective abelian sheaves are
$\mathcal{U}$-flabby.}

{\bf Corollary 2.} {\itshape For an abelian sheaf $F$ over $X$,
the following are equivalent.

$\mathrm{i)}$. $F$ is $\mathcal{U}$-flabby.

$\mathrm{ii)}$. For all $q>0$ and $U\in \mathcal{U}$, we have
$\mathcal{H}^{q}(F)(U)=0$ and therefore $H^{q}(F)(U)=0$.}

\subsection{Dagger formal schemes}
Let $A$ be a flat w.c.f.g. algebra over $R$. In this subsection,
we define $\mathrm{Spf}^{\dagger}(A)$.

$\mathrm{Spf}^{\dagger}(A)$ is a topology space together with a
structure sheaf $\mathcal{O}^{\dagger}$. The topology space is the
underlying space $X$ of $\mathrm{Spec}(A\otimes_{R}k)$.

For any $p\in X$, let us define
$\mathcal{O}^{\dagger}_{p}=\lim\limits_{\overrightarrow{\overline{f}\notin
p}}A\langle f^{-1}\rangle^{\dagger}$. Since $A\langle
f^{-1}\rangle^{\dagger}$ is flat over $A$,
$\mathcal{O}^{\dagger}_{p}$ is flat over $A$, too. Let us define
$\mathcal{O}^{\dagger}(U)$ to be the set of functions $s:
U\rightarrow\coprod\limits_{p\in
U}\mathcal{O}^{\dagger}_{p}\;\:(s(p)\in
\mathcal{O}^{\dagger}_{p})$ with the property that for each $p\in
U$, there is a $f\in A\;(\overline{f} \notin p)$ such that for any
$q\in X_{\overline{f}}\cap U$, $s(q)$ is the image of an element
of $A \langle f^{-1}\rangle^{\dagger}$ in
$\mathcal{O}^{\dagger}_{q}$.

It is obvious that $\mathcal{O}^{\dagger}$ is a sheaf.

{\bf Proposition 7.} {\itshape For any $f$ in $A$, we have
$\mathcal{O}^{\dagger}(X_{\overline{f}})=A\langle
f^{-1}\rangle^{\dagger}$.}

{\bf Proof.} We can define a homomorphism $\rho: A\langle
f^{-1}\rangle^{\dagger}\rightarrow\mathcal{O}^{\dagger}(X_{\overline{f}})$
in the natural way.

(1) At first, we show $\rho$ is injective. If $g,h\in A\langle
f^{-1}\rangle^{\dagger}$ such that $\rho(g)=\rho(h)$, then for any
$p\in X_{\overline{f}}$, $g$ and $h$ have the same image in
$\mathcal{O}^{\dagger}_{p}$. Then by the definition of
$\mathcal{O}^{\dagger}_{p}$, there is a $f_{p}\in
(f)\:(\overline{f}_{p}\notin p)$ such that $g$ and $h$ have the
same image in $A\langle f_{p}^{-1}\rangle^{\dagger}$. We choose a
finite number $\{f_{i}\}_{i=1} ^{r}$ of $\{f_{p}\}$ such that
$\{X_{\overline{f}_{i}}\}$ is a covering of $X_{\overline{f}}$.
Then from proposition 5, we obtain \begin{equation}0\rightarrow
A\langle f^{-1}\rangle^{\dagger}\rightarrow
\prod\limits_{i}A\langle
f_{i}^{-1}\rangle^{\dagger}\rightrightarrows\prod\limits_{i,j}A\langle
(f_{i}f_{j})^{-1}\rangle^{\dagger}.\label{*}\end{equation} So we
have $g=h$.

(2) Now we prove $\rho$ is surjective. Let $s$ be a section of $
\mathcal{O}^{\dagger}(X_{\overline{f}})$. By definition, for any
point $p\in X_{\overline{f}}$, there is a $f_{p}\in (f)
\;(\overline{f}_{p}\notin p)$ and $g\in A\langle
f_{p}^{-1}\rangle^{\dagger}$ such that for any $q\in
X_{\overline{f}_{p}}$, $s(q)$ is the image of $g$ in
$\mathcal{O}^{\dagger}_{q}$.

We can choose a finite subset $\{ f_{i}\}$ of $\{f_{p}\}$ such
that $\{X_{\overline{f}_{i}}\}$ is a covering of
$X_{\overline{f}}$. For each $f_{i}$, there is a $g_{i}\in
A\langle f_{i}^{-1}\rangle^{\dagger}$ such that $s(q)$ is the
image of $g_{i}$ in $\mathcal{O}^{\dagger}_{q}$ for each $q\in
X_{\overline{f}_{i}}$.

Since $\rho|_{A \langle (f_{i}f_{j})^{-1} \rangle^{\dagger} }$ is
injective, $g_{i}$ and $g_{j}$ have the same image in
$A\langle(f_{i}f_{j})^{-1}\rangle^{\dagger}$. Then from formula
(\ref{*}), we know $s\in \mathrm{Im}(\rho)$.\qed

{\bf Corollary 3.} {\itshape For a w.c.f.g. algebra $A$ over $R$,
we have $\mathcal{O}^{\dagger}(\mathrm{Spf}^{\dagger}(A))=A$.}

We call such a pair $(\mathrm{Spec}(A\otimes_{R}k ),
\mathcal{O}^{\dagger})$ an affine dagger formal scheme
$\mathrm{Spf}^{\dagger}(A)$.

{\bf Corollary 4.} {\itshape Let $U\subset
X=\mathrm{Spf}^{\dagger}(A)$ be an open subset with
$\{X_{\overline{f}_{i}}\}$ a covering of $U$. Then we have the
following exact sequence}
$$0\rightarrow \mathcal{O}^{\dagger}(U) \rightarrow \prod_{i}A\langle f_{i}^{-1}\rangle^{\dagger}\rightrightarrows\prod_{i,j}A\langle(f_{i}f_{j})^{-1}\rangle^{\dagger}.$$

{\bf Proposition 8.} {\itshape Affine dagger formal schemes are
locally ringed spaces.}

{\bf Proof.} Let $X=\mathrm{Spf}^{\dagger}(A)$. From Proposition
7, we know for each $p\in X$, the stalk of
$\mathcal{O}^{\dagger}_{X}$ at $p$ is $\mathcal{O}^{\dagger}_{p}$.
We need to prove $\mathcal{O}^{\dagger}_{p}$ is a local ring.

Since $\pi\in \mathrm{rad}(A\langle f_{i}^{-1}\rangle^{\dagger})$,
and $\mathcal{O}^{\dagger}_{p}$ is the limit of $A\langle
f^{-1}\rangle^{\dagger}$ ($\bar{f} \notin p$), we see $\pi\in
\mathrm{rad}(\mathcal{O}^{\dagger}_{p})$. Since
$\mathcal{O}^{\dagger}_{p}/\pi \mathcal{O}^{\dagger}_{p}$ is the
stalk of
$\mathcal{O}_{\bar{X}}=\mathcal{O}^{\dagger}_{X}\otimes_{R}k$ at
$p$, it is a local ring. Therefore, $\mathcal{O}^{\dagger}_{p}$ is
also a local ring. \qed

Now we can define dagger formal schemes.

{\bf Definition.} {\itshape A dagger formal scheme is a locally
ringed space $(X,\mathcal{O}^{\dagger}_{X})$ in which every point
has an open neighborhood $U$, such that
$(U,\mathcal{O}^{\dagger}_{X}|U)$ is an affine dagger formal
scheme. A morphism of dagger formal schemes is a morphism as
locally ringed spaces.}

For a given $A$-module $M$, we can associate to it a sheaf
$M^{\triangle}$ of
$\mathcal{O}^{\dagger}_{\mathrm{Spf}^{\dagger}(A)}$-modules over
$\mathrm{Spf}^{\dagger}(A)$ (say simply a
$\mathcal{O}^{\dagger}_{\mathrm{Spf}^{\dagger}(A)}$-module).

For $U\subset \mathrm{Spf}^{\dagger}(A)$, define
$M^{\triangle}(U)$ to be the set of functions $s: U\rightarrow
\coprod\limits_{p\in U}\mathcal{O}^{\dagger}_{p}\otimes_{A} M\;\:
(s(p)\in \mathcal{O}^{\dagger}_{p}\otimes_{A} M)$ with the
property that for each $p\in U$, there is a $f\in A\;
(\overline{f}\notin p)$ such that for each $q\in
X_{\overline{f}}\cap U$, $s(q)$ is the image of an element of
$(A\langle f^{-1}\rangle^{\dagger})\otimes_{A}M$. By definition,
we know that $M^{\triangle}$ is a
$\mathcal{O}^{\dagger}_{\mathrm{Spf}^{\dagger}(A)}$-module.

{\bf Proposition 9.} {\itshape For any $f\in A$, we have}
$$M^{\triangle}(X_{\overline{f}})=A\langle f^{-1}\rangle^{\dagger}\otimes_{A}M.$$

Using the exact sequence given in Proposition 5, we can prove this
proposition in the same way as Proposition 7.

Let $\widetilde{\mathcal{O}}^{\dagger}$ be
$\mathcal{O}^{\dagger}\otimes_{R}K$, and
$\widetilde{M}^{\triangle}$ be $M^{\triangle}\otimes_{R}K$.

Now we can define coherent $\mathcal{O}^{\dagger}_{X}$-modules and
quasi-coherent $\mathcal{O}^{\dagger}_{X}$-modules as in scheme
case. We have the following proposition.

{\bf Proposition.} {\itshape Let $h: X\rightarrow Y$ be a morphism
of Noetherian dagger formal schemes over $R$. Let $\mathcal{G}$ be
a sheaf of $\mathcal{O}^{\dagger}_{Y}$-modules and $\mathcal{F}$
be $\mathcal{O}^{\dagger}_{X}$-modules.

$\mathrm{(i)}$. If $\mathcal{G}$ is quasi-coherent, then
$h^{*}(\mathcal{G})$ is also quasi-coherent.

$\mathrm{(ii)}$. If $\mathcal{G}$ is coherent, then
$h^{*}(\mathcal{G})$ is also coherent.

$\mathrm{(iii)}$. If $\mathcal{F}$ is quasi-coherent, then
$h_{*}(\mathcal{F})$ is also quasi-coherent.}

Let $A$ be a w.c.f.g. algebra over $R$ with $f\in A$ and
$f_{1},...,f_{r}\in (f)$ satisfy formula $(*)$ given in section 2.
Let $U=X_{\overline{f}}$ and $U_{i}=X_{\overline{f}_{i}}$. Then
Tate's acyclic theorem says
$$0\rightarrow A\langle f^{-1}\rangle^{\dagger}\otimes_{R}K\rightarrow\prod_{i}A\langle f_{i}^{-1}\rangle^{\dagger}\otimes_{R}K\rightarrow \breve{C}^{1}(\{U_{i}\rightarrow U\},\widetilde{\mathcal{O}}^{\dagger})\rightarrow\cdot\cdot\cdot\rightarrow 0.$$
Since $\breve{C}^{i}(\{U_{i}\rightarrow
U\},\widetilde{\mathcal{O}}^{\dagger})$ are all flat $A$-modules,
we get{\scriptsize
$$0\rightarrow (A\langle f^{-1}\rangle^{\dagger}\otimes_{R}K)\otimes_{A}M\rightarrow\prod_{i}(A\langle f_{i}^{-1}\rangle^{\dagger}\otimes_{R}K)\otimes_{A}M\rightarrow \breve{C}^{1}(\{U_{i}\rightarrow U\},\widetilde{M}^{\triangle})\rightarrow\cdot\cdot\cdot\rightarrow 0.$$
}Therefore we obtained the following proposition.

{\bf Proposition 10.} {\itshape For $j>0$, we have
$\breve{H}^{j}(\{U_{i}\rightarrow
U\},\widetilde{M}^{\triangle})=0$.}

As a consequence of this proposition and Corollary 2, we get the
following proposition.

{\bf Proposition 11.} {\itshape Let $X=\mathrm{Spf}^{\dagger}(A)$.
For $i>0$, we have
$H^{i}(X_{\overline{f}},\widetilde{M}^{\triangle})=0$.}


\subsection{Morphism}
In this subsection, we study morphisms of dagger formal schemes.

{\bf Lemma 3.} {\itshape Every homomorphism of w.c.f.g
$R$-algebras is continuous.}

{\bf Proof.} Let $\varphi:A_{1}\rightarrow A_{2}$ be a
homomorphism of w.c.f.g $R$-algebras. From
$\varphi(\pi^{i}A_{1})\subset \pi^{i}A_{2}$, we see $\varphi$ is
continuous.\qed

{\bf Proposition 12.} {\itshape Let $A,B$ be two w.c.f.g.
$R$-algebras.

$\mathrm{(i)}$. If $\rho:A\rightarrow B$ is a homomorphism of
$R$-algebras, then $\rho$ induces a natural morphism of locally
ringed spaces
$$( \varphi, \varphi ^{\sharp} ):(\mathrm{Spf}^{\dagger}(B),\mathcal{O}_{B}^{\dagger})\rightarrow (\mathrm{Spf}^{\dagger}(A),\mathcal{O}_{A}^{\dagger})$$

$\mathrm{(ii)}$. Every morphism of locally ringed spaces from
$\mathrm{Spf}^{\dagger}(B)$ to $\mathrm{Spf}^{\dagger}(A)$ is
induced by a $R$-morphism $\rho: A\rightarrow B$as in
$\mathrm{(i)}$.}

{\bf Proof.} (i) is easy. We only need to prove (ii).

Let ($\varphi,\varphi^{\sharp}$) be a morphism of locally ringed
spaces from $Y=\mathrm{Spf}^{\dagger}(B)$ to
$X=\mathrm{Spf}^{\dagger}(A)$. It induces a morphism
($\bar{\varphi},\bar{\varphi}^{\sharp}$) of locally ringed spaces
from $\bar{Y}=\mathrm{Spec}(B\otimes _{R}k)$ to
$\bar{X}=\mathrm{Spec}(A\otimes _{R}k)$. It is well know that
($\varphi,\bar{\varphi}^{\sharp}$) is induced from a homomorphism
$\bar{\rho}$ of $k$-algebras $A\otimes_{R} k\rightarrow
B\otimes_{R}k$. Therefore, we see for any $f\in A$
\begin{equation}\varphi(Y_{\bar{f}})=X_{\overline{\rho(f)}}.\label{easy}\end{equation}

Let $\rho$ be $\varphi^{\sharp}(Y): A\rightarrow B$. From
(\ref{easy}), we get the following commutative diagram
\[
\begin{CD}
A @>\rho>> B \\
@VVV @VVV\\
A\langle f^{-1}\rangle^{\dagger}@>\varphi^{\sharp}(X_{\bar{f}})>>
B\langle\rho(f)^{-1}\rangle^{\dagger}.
\end{CD}
\]

From lemma 3, we know $\rho_{f}$ and
$\varphi^{\sharp}(X_{\bar(f)})$ are continuous. Since they
restrict to the dense subalgebra $A_{f}$ are the same, they are
the same. \qed

Let $k$ be a finite field with $q$ elements. Assume
$R=\mathrm{W}(k)$. There are Frobenius actions on $k$-algebras
$A\rightarrow A$ $(x\mapsto x^{q})$, and on $k$-schemes.

A $R$-morphism $\mathrm{F}$ from a dagger formal scheme $X$ to
itself is called a Frobenius of $X$ if it induces the Frobenius
over $\bar{X}=X\otimes_{R}k$.

By definition, a Frobenius over $X$ is equivalent to a collection
$\{\mathrm{F}_{U}\}$ of homomorphisms
$\mathrm{F}_{U}:\mathcal{O}^{\dagger}_{X}(U)\rightarrow
\mathcal{O}^{\dagger}_{X}(U)$ for each $U \subset X$ such that
$\mathrm{F}_{U}\otimes_{R}k$ is the Frobenius of
$\mathcal{O}^{\dagger}_{X}(U)\otimes_{R}k$ and the following
diagram is commutative $(U\subset V)$
\[
\begin{CD}
\mathcal{O}_{X}^{\dagger}(V) @>\mathrm{F}_{V}>>
\mathcal{O}_{X}^{\dagger}(V)\\
@V\mathrm{res}VV@V\mathrm{res}VV\\
\mathcal{O}_{X}^{\dagger}(U) @>\mathrm{F}_{U}>>
\mathcal{O}_{X}^{\dagger}(U).
\end{CD}
\]

\subsection{Product}
Let $A$ and $B$ be two flat w.c.f.g. algebras over $R$. In this
subsection, we define $A\hat{\otimes}^{\dagger}_{R}B$ which is
also a w.c.f.g. algebra.

{\bf Lemma 4.} {\itshape Let $C$ be a finite generated $R$-algebra
with $\breve{C}$ its $\pi$-adic completion. Let $D$ be a c.f.g.
algebra. Then any homomorphism $\varphi: C\rightarrow D$ extends
to a unique homomorphism $\varphi^{\mathrm{v}}:
\breve{C}\rightarrow D$. Moreover when $\varphi$ is injective,
$\varphi^{\mathrm{v}}$ is also injective.}

{\bf Proof.} Existence of $\varphi^{\mathrm{v}}$ comes from
$\breve{C}=\varprojlim C/\pi^{i} C$ and $D=\varprojlim D/\pi^{i}
D$. For uniqueness of $\varphi^{\mathrm{v}}$, one only needs to
note that a homomorphism of two c.f.g. $R$-algebras is
continuous.\qed

We know $A=\varinjlim A_{n}^{\mathrm{v}}$ with
$A_{n}^{\mathrm{v}}$ c.f.g. algebras over $R$ such that
$A_{n}^{\mathrm{v}}\rightarrow A_{n+1}^{\mathrm{v}}$ is injective.
Since torsion free $R$-modules are flat over $R$,
$A_{n}^{\mathrm{v}}$ are all flat over $R$. The same holds for
$B$.

Then we define $A\hat{\otimes}^{\dagger}_{R}B=\varinjlim
(A_{n}^{\mathrm{v}}\hat{\otimes} B_{m}^{\mathrm{v}})$.

By Lemma 4, if $A=\varinjlim A_{n}^{\mathrm{v}}$ and $A=\varinjlim
 {A'}_{n}^{\mathrm{v}}$, then $A_{n}^{\mathrm{v}}$ is contained in
some ${A'}_{m'}^{\mathrm{v}}$ and $A_{m}^{\mathrm{v}}$ is
contained in some ${A'}_{n'}^{\mathrm{v}}$. Therefore, the
definition of $A\hat{\otimes}^{\dagger}_{R}B$ does not depend on
the choice of $\{A_{n}^{\mathrm{v}} \}$ and $\{ B_{m}^{\mathrm{v}}
\}$.

Assume $A=\mathrm{T}_{n}^{\dagger}/I_{1}$ and
$B=\mathrm{T}_{m}^{\dagger}/I_{2}$. Let $i_{n}:
\mathrm{T}_{n}^{\dagger}\rightarrow\mathrm{T}_{n+m}^{\dagger}$ be
the homomorphism maps variants of $\mathrm{T}_{n}^{\dagger}$ to
first $n$ variants, and let
$j_{m}:\mathrm{T}_{m}^{\dagger}\rightarrow\mathrm{T}_{n+m}^{\dagger}$
be the homomorphism maps variants of $\mathrm{T}_{n}^{\dagger}$ to
last $m$ variants. Then
$A\hat{\otimes}^{\dagger}_{R}B=\mathrm{T}_{n+m}^{\dagger}/I$ with
$I$ the ideal of $\mathrm{T}_{n+m}^{\dagger}$ generated by
$i_{n}(I_{1})$ and $j_{m}(I_{2})$.

Let $i_{A}:A\rightarrow A\hat{\otimes}^{\dagger}_{R}B$ and
$j_{B}:B\rightarrow A\hat{\otimes}^{\dagger}_{R}B$ be two
canonical inclusion homomorphisms. Let $C$ be a w.c.f.g. algebra
over $R$.

{\bf Lemma 5.} {\itshape If $\varphi_{1}: A\rightarrow C$ and
$\varphi_{2}:B\rightarrow C$ are two homomorphisms of
$R$-algebras, then there is a unique homomorphism
$\varphi:A\hat{\otimes}^{\dagger}_{R}B\rightarrow C $ such that
$\varphi\circ i_{A}=\varphi _{1}$ and $\varphi\circ j_{B}=\varphi
_{2}$.}

{\bf Proof.} Assume that $A=\varinjlim A_{n}^{\mathrm{v}}$,
$B=\varinjlim B_{m}^{\mathrm{v}}$ and $C=\varinjlim
C_{l}^{\mathrm{v}}$ as in section 3.4. Then there is a $l$, such
that $\varphi_{1}(A_{n}^{\mathrm{v}})$ and
$\varphi_{2}(B_{m}^{\mathrm{v}})$ are contained in
$C_{l}^{\mathrm{v}}$. By Lemma 4, there is a $\varphi_{n,m}:
A_{n}^{\mathrm{v}}\hat{\otimes}B_{m}^{\mathrm{v}}\rightarrow
C_{l}^{\mathrm{v}}$, such that $\varphi_{n,m}\circ i_{A}=\varphi
_{1}|_{A_{n}}$ and $\varphi_{n,m}\circ j_{B}=\varphi
_{2}|_{B_{m}}$. Taking limit, we get what we desire. \qed

For $X=\mathrm{Spf}^{\dagger}(A)$ and
$Y=\mathrm{Spf}^{\dagger}(B)$, we can define $X\times
Y=\mathrm{Spf}^{\dagger}(A\hat{\otimes}^{\dagger}_{R}B)$. By
glueing, we can define $X\times Y$ for separated flat dagger
formal schemes $X$ and $Y$.

\section{De Rham cohomology and Lefchetz fixed points theorem}
\subsection{Differential modules and de Rham cohomology}
In this subsection, we recall the concept of differentaaial
modules.

Let $A$ be a w.c.f.g. algebra over $R$.

The module of differential forms of $A$ over $R$ is a finite
$A$-module $\Omega^{1}_{A/R}$ together with a $R$-derivation
$d_{A/R}$, which is universal in the following sense: for any
finite $A$-module $M$, the canonical map
$$\mathrm{Hom}_{A}(\Omega^{1}_{A/R},M)\xrightarrow{\mbox{~}}\mathrm{Der}_{R}(A,M),\;\;\varphi\mapsto \varphi\circ d_{A/R}$$
is bijective.

Let $m:A\hat{\otimes}^{\dagger}A\rightarrow A$ be the ``diagonal
homomorphism" defined by $m(b\otimes b')=bb'$. By Lemma 5, this
homomorphism makes sense. Let $I$ be the kernel of $m$. Then
$I/I^{2}$ inherits a structure of $A$-module. Define a map $d:
A\rightarrow I/I^{2}$ by $db=1\otimes b-b\otimes 1
\:(\mathrm{mod}I^{2})$.

{\bf Proposition 13.} {\itshape $(I/I^{2}, d)$ is
$(\Omega^{1}_{A/R},d_{A/R})$.}

{\bf Proof.} (1). Let $M$ be a finite $A$-module. Define $A*M$ by
$$(a_{1},m_{1})+(a_{2},m_{2})=(a_{1}+a_{2},m_{1}+m_{2})$$
and
$$(a_{1},m_{1})\cdot (a_{2},m_{2})=(a_{1}a_{2},a_{1}m_{2}+a_{2}m_{1})$$
for $a_{1},a_{2}\in A$ and $m_{1},m_{2}\in M$. Then $A*M$ is a
w.c.f.g. algebra over $R$.

It is sufficient to prove this fact in the case
$A=\mathrm{T}_{n}^{\dagger}$ and $M=\mathrm{T}_{n}^{\dagger\oplus
r}$. In this case, $A*M=\mathrm{T}^{\dagger}_{n+r}/\sim$, where
$\mathrm{T}^{\dagger}_{n+r}=R\langle\xi_{1},...,\xi_{n},\xi_{n+1},...,\xi_{n+r}\rangle^{\dagger}$,
and $\sim$ is generated by $\{\:\xi_{i}\xi_{j}\:|i\geq n+1, j\geq
n+1\}$.

(2). We need to prove the fact that if $D$ is an $R$-derivation of
$A$ into an $A$-module $M$, then there is a unique $A$-linear map
$f:I/I^{2}\rightarrow M$ such that $D=fd$.

(3). Let $I'$ be $\mathrm{Ker}(A\otimes A\rightarrow A)$, then
$I=I'(A\hat{\otimes}^{\dagger}A)$. In the following, we prove this
assertion.

Since $A=\varinjlim A_{n}^{\mathrm{v}}$, $I$ is generated by
elements in
$I_{n}=\mathrm{Ker}(\hat{A}_{n}\hat{\otimes}\hat{A}_{n}\rightarrow\hat{A}_{n})$
for $n$ large enough. Let
$I'_{n}=\mathrm{Ker}(\hat{A}_{n}\otimes\hat{A}_{n}\rightarrow\hat{A}_{n})$.
To show $I=I'(A\hat{\otimes}^{\dagger}A)$, it is sufficient to
show $I_{n}=I'_{n}(\hat{A}_{n}\hat{\otimes}\hat{A}_{n})$.

Since $\hat{A}_{n}\otimes\hat{A}_{n}$ is dense in
$\hat{A}_{n}\hat{\otimes}\hat{A}_{n}$,
$\hat{A}_{n}=\hat{A}_{n}\otimes\hat{A}_{n}/I'_{n}$ is dense in
$\hat{A}_{n}\hat{\otimes}\hat{A}_{n}/I'_{n}(\hat{A}_{n}\hat{\otimes}\hat{A}_{n})$.
Since both of them are c.f.g. algebras, them must be the same.
Therefore, $I_{n}=I'_{n}(\hat{A}_{n}\hat{\otimes}\hat{A}_{n})$.

(4). By (3), we know $I/I^{2}$ is generated by $\{dy|y\in A\}$ as
$A$-module. Then we get the uniqueness of $f$.

(5). By Lemma 5, from the following two homomorphisms
$$\phi_{1}: A\rightarrow A*M,\;\;\: \phi_{1}(x)=(x,0)$$
and
$$\phi_{2}: A\rightarrow A*M,\;\;\: \phi_{2}(x)=(x, D(x)),$$
we get a homomorphism
$$\phi: A\hat{\otimes}^{\dagger}A\rightarrow A*M$$
whose restriction on $A\otimes A$ is
$$\phi(x\otimes y)=(xy,xD(y)).$$
Then $\phi(I')\subseteq M$, and so $\phi(I)\subseteq M$. Since
$M^{2}=0$, $\phi(I^{2})=0$. So $\phi$ induces
$\bar{\phi}:(A\hat{\otimes}^{\dagger}A)/I^{2}=A*(I/I^{2})\rightarrow
A*M$ which maps $dy$ to $(0, D(y))$.

Thus the restriction of $\bar{\phi}$ on $I/I^{2}$ gives an
$A$-linear map $f:I/I^{2}\rightarrow M$ such that $f\circ d=D$.
\qed

If
$A=\mathrm{T}_{n}^{\dagger}=R\langle\xi_{1},...,\xi_{n}\rangle$,
then
$\Omega^{1}_{A/R}=\mathrm{T}_{n}^{\dagger}d\xi_{1}\oplus\cdot\cdot\cdot\oplus\mathrm{T}_{n}^{\dagger}
d \xi_{n}$. If $A=\mathrm{T}_{n}^{\dagger}/I$, then
$\Omega^{1}_{A/R}=\Omega^{1}_{\mathrm{T}_{n}^{\dagger}/R}/\sim$
with $\sim$ generated by
$I\Omega^{1}_{\mathrm{T}_{n}^{\dagger}/R}$ and
$\mathrm{T}_{n}^{\dagger} df$ ($f\in I$).

For a dagger formal scheme $X$ over $R$, define
$$\Omega^{1}_{X/R}:=\Delta^{*}(\mathcal{I}/\mathcal{I}^{2})$$
where $\Delta:X\rightarrow X\times X$ is the diagonal map and
$\mathcal{I}$ the ideal defining the diagonal map.

From Proposition 13, we know when $X=\mathrm{Spf}^{\dagger}(A)$,
$\Omega^{1}_{X/R}$ is the sheaf associated to the $A$-module
$\Omega^{1}_{A/R}$.

We have the following proposition whose proof is the same as
\cite{Put}(2.3).

{\bf Proposition 14.} {\itshape Let $X$ be a flat and separated
dagger formal scheme. If $\bar{X}=X\otimes_{R}k$ is regular, then
$\Omega^{1}_{X/R}$ is a locally free
$\mathcal{O}_{X}^{\dagger}$-module with rank
$\mathrm{dim}(\bar{X})$.}

When $X$ satisfies the condition of Proposition 14, we define
$\Omega^{i}_{X/R}=\wedge^{i} \Omega^{1}_{X/R}$, and get the de
Rham complexs $(\Omega^{i}_{X/R},d)$ and
$(\widetilde{\Omega^{i}_{X/R}},\tilde{d})$ with
$\widetilde{\Omega^{i}_{X/R}}=\Omega^{i}_{X/R}\otimes_{R}K$. We
write $d$ in place of $\tilde{d}$ for simple. Then we define
$H_{dR}^{i}(X,K):=\mathbb{H}^{i}(X,(\widetilde{\Omega^{i}_{X/R}},d))$.
From Proposition 11, we see when $X=\mathrm{Spf}^{\dagger}(A)$,
$H_{dR}^{i}(X,K)$ is exactly the Monsky-Washnitzer cohomology
$H_{\small\mathrm{MW}}^{i}$ defined to be
$H^{i}(\Omega^{i}_{A/R}\otimes_{R}K,d))$.

From now on, assume $R=\mathrm{W}(k)$.

When $X$ is a flat and separated dagger formal scheme over $R$
with a Frobenius $\mathrm{F}$ such that $\bar{X}$ is regular and
integral, then $\mathrm{F}$ induces $\mathrm{F}^{*}$ over
$H_{dR}^{i}(X,K)$.

\subsection{nuclear operator and operator $\psi$}
\subsubsection{nuclear operator}
{\bf Definition.} {\itshape A $K$-linear map $L:M\rightarrow M$ is
called nuclear, if the following two conditions hold.\\
$(\mathrm{i})$. For every $\lambda\neq 0$ in $K^{ac}$ (the
algebraic closure of $K$) with $g$ the minimal polynomial of
$\lambda$ over $K$, $\cup(\mathrm{Ker}(g(L)^{m}))$ is of finite
dimension.\\
$(\mathrm{ii})$. The nonzero eigenvalues of $L$, form a finite set
or a sequence with a limit 0.}

From (i), we see $M=V\oplus W$ with $V$, $W$ vector spaces
invariant under $L$ such that $W=\cup(\mathrm{Ker}(g(L)^{m}))$ and
$g(L)$ is bijection over $V=\cap \mathrm{Im}(g(L)^{m})$.

We can define trace $\mathrm{tr}(L)$ for nuclear operator $L$. Let
$M_{l}$ be the sum of the generalized eigenspaces of $L$ with
eigenvalues $\lambda$ ($|\lambda|\geq |\pi|^{l}$). Then
$\mathrm{dim}M_{l}<\infty$. Define
$\mathrm{tr}(L^{s})=\lim\limits_{l\rightarrow\infty}\mathrm{tr}(L^{s}|_{M_{l}})$
for positive integer $s$, and
$\mathrm{det}(1-tL)=\lim\limits_{l\rightarrow\infty}\mathrm{det}(1-tL|_{M_{l}})$.

For nuclear operators, we have the following two lemmas.

 {\bf Lemma 6.} {\itshape Let $L_{i}:
M_{i}\rightarrow M_{i}$ $(i=1,2)$ be nuclear. Assume linear map
$\alpha:M_{1}\rightarrow M_{2}$ satisfy $\alpha
L_{1}=L_{2}\alpha$. Then the induced maps $L_{0}$ on $\mathrm{Ker}
(\alpha)$ and $L_{3}$ on $\mathrm{Coker} (\alpha)$ are nuclear.
Moreover,
$$\prod_{i=0}^{3}\mathrm{det}(1-tL_{i})^{(-1)^{i}}=1,$$
and $$\sum_{i=0}^{3}\mathrm{tr}(L_{i}^{s})=0,$$ where $s$ is a
positive integer.}

{\bf Lemma 7.} {\itshape Let $L$ be a $K$-linear map over $M$.
Assume $M_{1}$ is a $K$-linear subspace of $M$ fixed by $L$. Then
$L$ induces $L_{1}$ on $M_{1}$ and $L_{2}$ on $M/M_{1}$. Then $L$
is nuclear if and only if both $L_{1}$ and $L_{2}$ are nuclear. If
so, then we have the following formulas
$$\mathrm{det}(1-tL)=\mathrm{det}(1-tL_{1})\cdot\mathrm{det}(1-tL_{2}),$$
and
$$\mathrm{tr}(L^{s})=\mathrm{tr}(L_{1}^{s})+\mathrm{tr}(L_{2}^{s}),$$
where $s$ is a positive integer.}

\subsubsection{Operator $\psi$}

{\bf Proposition.} (\cite{Put}) {\itshape Let $B\subset A$ denote
a finite ring extension of w.c.f.g. algebras. Suppose both
$\bar{A}$ and $\bar{B}$ are regular, and both $A$ and $B$ are
integral and flat over $R$. Then there exists a ``trace map"
$tr_{A/B}:\Omega^{i}(A)\rightarrow \Omega^{i}(B)$.}

$Tr_{B/A}$ is defined by
$$
\Omega^{i}(A) \rightarrow \Omega^{i}(A)\otimes_{A}\mathrm{Qt}(A)
\xrightarrow[\cong]{}  \Omega^{i}(B)\otimes_{A}\mathrm{Qt} (A)
\xrightarrow[ \mathrm{id}\otimes tr]{}
\Omega^{i}(B)\otimes_{B}\mathrm{Qt}(B).
$$
In \cite{Put}, one shows $Tr_{A/B}$ maps $\Omega_{i}(A)$ into
$\Omega_{i}(B)$. Moreover, one can show $tr_{A/B}$ commutes with
$d$.

Let $A$ be a w.c.f.g. algebra over $R$ satisfying the above
proposition with $\bar{A}$ a regular and integral algebra of
dimension $n$. Let $\mathrm{F}$ be a Frobenius over $A$. Since
$[A:\mathrm{F}(A)]=[A\langle
f^{-1}\rangle^{\dagger}:\mathrm{F}(A\langle
f^{-1}\rangle^{\dagger})]=q^{n},$ we have the following
commutative diagram
\begin{equation}
\begin{CD}
\Omega^{i}(A)\otimes_{R}K @>tr>> \Omega^{i}(\mathrm{F}(A))\otimes_{R}K\\
@V res VV @V res VV\\
\Omega^{i}(A\langle f^{-1}\rangle^{\dagger})\otimes_{R}K @>tr>>
\Omega^{i}(\mathrm{F}(A\langle
f^{-1}\rangle^{\dagger}))\otimes_{R}K. \label{define trace}
\end{CD}
\end{equation}

Let $X$ be a flat and integral dagger formal scheme with a
Frobenius $\mathrm{F}$ and regular and integral reduction
$\bar{X}$. From commutative diagram (\ref{define trace}), we get a
trace map $tr_{X}: X=(\bar{X},
\mathcal{O}_{X}^{\dagger})\rightarrow
X^{q}=(\bar{X},\mathrm{F}(\mathcal{O}_{X}^{\dagger}))$, with
$\mathrm{F}(\mathcal{O}_{X}^{\dagger})(U)=\mathrm{F}_{U}(\mathcal{O}_{X}^{\dagger}(U))$
($U\subset X$).

We define operator $\psi=\mathrm{F}^{-1}\circ tr_{X}$:
$$\Omega^{i}(X,K)\xrightarrow{tr}\Omega^{i}(X^{(q)},K)\xleftarrow[\sim]{\mbox{F}}\Omega^{i}(X,K).$$

{\bf Proposition 15.} (\cite{Put}) {\itshape

$\mathrm{(i)}$. For any open subset $U\subset X$, $a\in
\mathcal{O}_{X}^{\dagger}(U),\omega \in \Omega^{i}(U,K)$, we have
$\psi ( \mathrm{F}(a) \omega) =a\psi ( \omega ).$

$\mathrm{(ii)}$. $\psi$ commutes with $d$.

$\mathrm{(iii)}$. $\psi\circ F=q^{n}$.}

From \cite{Put}, we know $\psi_{U}$ is nuclear for any open subset
$U$ of $X$. Let $\psi^{*}$ be the linear morphisms on
$H_{dR}^{i}(X,K)$ induced by $\psi$.

{\bf Proposition 16.} (\cite{Put}) {\itshape Let $A$ be a flat
w.c.f.g. algebra over $R$ with $\bar{A}=A\otimes_{R}k$ regular and
integral. Then $\mathrm{F}^{*}$ is  bijection over $H_{\small
\mathrm{MW}}(A,K)$ and $\psi^{*}=q^{n}(\mathrm{F}^{*})^{-1}$ is
nuclear.}

\subsection{Lefchetz fixed points theorem}
In \cite{Put} (4.1), the following result is stated.

{\bf Theorem 2.} {\itshape Let $A$ be a flat w.c.f.g. algebra over
$R$ with $\bar{A}=A\otimes_{R}k$ regular and integral of dimension
$n$. Let $N(\bar{A})$ denote the number of $k$-homomorphisms
$\bar{A}\rightarrow k=\mathbb{F}_{q}$. Then}
\begin{equation}
N(\bar{A})=\sum_{i=0}^{n}(-1)^{i}\mathrm{tr}(\psi^{*}|H^{i}_{\small
\mathrm{MW}}(A,K)).
\end{equation}

In this subsection, we show the following generalization of
Theorem 2.

{\bf Theorem 3.} {\itshape Let $X$ be a flat separated and
Noetherian dagger formal scheme with $\bar{X}=X\otimes_{R}k$
regular and integral of dimension $n$. Let $N(\bar{X})$ denote the
number of $k$-points of $\bar{X}$. Then $\psi^{*}$ is nuclear over
$H_{dR}^{i}(X,K)$ and the following formula  holds}
\begin{equation}
N(\bar{X})=\sum (-1)^{i}\mathrm{tr}(\psi^{*}|H^{i}_{dR}(X,K)).
\label{Lefschetz}
\end{equation}

When $X=\mathrm{Spf}^{\dagger}(A)$, $H_{dR}^{i}(X,K)=H_{\small
\mathrm{MW}}(\bar{A},K)$. So Theorem 2 is a special case of
Theorem 3.

{\bf Proof.} (1). Take an affine covering
$\mathcal{U}=\{U_{i}\}_{i=1}^{r}$ of $X$. Since $X$ is separated,
intersections of $\{U_{i}\}'$s are again affine.

Then from Proposition 11 and EGA III (12.4.7), we obtain
\begin{equation}
\breve{\mathbb{H}}^{\cdot}(\mathcal{U},\widetilde{\Omega}^{\cdot})\xrightarrow[\sim]{}\mathbb{H}^{\cdot}(X,\widetilde{\Omega}^{\cdot})
\end{equation}
So we can calculate
$\mathbb{H}^{\cdot}(X,\widetilde{\Omega}^{\cdot})$ by spectral
sequence.

Write $K^{p,q}=\breve{C}^{p}(\mathcal{U},\widetilde{\Omega}^{q})$.
Let
\[
F^{p}K^{i,j}=\left \{
\begin{array}{l@{\quad\quad}l}
K^{i,j} & \mathrm{if}\; i \geq p \\
0 & \mathrm{if}\;i\leq p
\end{array}
\right.
\]

Let $sK^{l}=\oplus_{p+q=l}K^{p,q}$. We define $sF^{p}H^{\cdot}$ in
a similar way.

Let
$E_{1}^{p,q}=\breve{C}^{p}(\mathcal{U},H^{q}(\cdot,\widetilde{\Omega}^{\cdot}))=\breve{C}^{p}(\mathcal{U},H^{q}_{\small
\mathrm{MW}}(\cdot\:))$, then
\begin{equation}
E_{1}^{p,q}\Longrightarrow
\breve{\mathbb{H}}^{\cdot}(\mathcal{U},\widetilde{\Omega}^{\cdot})\label{spectral
sequence}
\end{equation}
with
$F^{p}\breve{\mathbb{H}}^{l}(\mathcal{U},\widetilde{\Omega}^{\cdot})=\mathrm{Im}(H^{l}(sF^{p}K)\rightarrow
H^{l}(sK))$.

(2). We show $\psi^{*}$ is nuclear.

From Proposition 16, we see $\psi$ is nuclear on each term of
$E_{1}^{p,q}$. Then by Lemma 6 and Lemma 7, we obtain that
$\psi^{*}$ is nuclear on each term of $E_{m}^{p,q}$ with $m$ a
positive integer. Therefore, from (\ref{spectral sequence}), we
see that $\psi^{*}$ is nuclear on $H_{dR}^{n}(X,K)$.

(3). From
$$N(\bar{X})=\sum\limits_{i_{0}<...<i_{j}}(-1)^{j}N(\bar{U}_{i_{1}...i_{j}})$$
and
$$N(\bar{U}_{i_{0}...i_{j}})=\sum (-1)^{i}\mathrm{tr}(\psi^{*}|H^{i}_{\small \mathrm{MW}}(U_{i_{0}...i_{j}},K)),$$
we get
$$N(\bar{X})=\sum_{p,q}(-1)^{p+q}\mathrm{tr}(\psi^{*}|E_{1}^{p,q}).$$
From Lemma 6 and Lemma 7, we get
$$N(\bar{X})=\sum_{p,q}(-1)^{p+q}\mathrm{tr}(\psi^{*}|E_{m}^{p,q}),$$
for all positive integer $m$. Then from Lemma 6, Lemma 7 and
(\ref{spectral sequence}), We get (\ref{Lefschetz}). \qed

From Proposition 16 and (\ref{spectral sequence}), we get the
following lemma easily.

{\bf Lemma 8.} {\itshape $\mathrm{F}^{*}$ is bijection on
$H_{dR}^{i}(X,K)$ and $\psi^{*}=q^{n}(\mathrm{F}^{*})^{-1}$.}

{\bf Theorem 4.} {\itshape Let $X$ be a flat separated and
Noetherian dagger formal scheme with $\bar{X}=X\otimes_{R}k$
regular and integral of dimension $n$. Let $N_{s}(\bar{X})$ denote
the number of $\;\mathbb{F}_{q^{s}}$-points of $\bar{X}$. Then the
following formula holds}
\begin{equation}
N_{s}(\bar{X})=\sum
(-1)^{i}\mathrm{tr}((q^{n}(\mathrm{F}^{*})^{-1})^{s}|H^{i}_{dR}(X,K)).
\end{equation}

\end{document}